\documentclass[10pt,a4paper,twoside]{article}

\usepackage{amssymb,amsmath,multirow}

\newtheorem{theorem}{Theorem}[section]

\newtheorem{corollary}[theorem]{Corollary}

\newtheorem{lemma}[theorem]{Lemma}

\newtheorem{remark}[theorem]{Remark}

\DeclareRobustCommand{\stirling}{\genfrac\{\}{0pt}{}}

    \makeatletter
    \let\@fnsymbol\@arabic
    \makeatother

\title{Some Characterizations of Exponential Distribution Based on Order Statistics}
\author{Bojana Milo\v sevi\' c\footnote{bojana@matf.bg.ac.rs}, Marko Obradovi\' c\footnote{marcone@matf.bg.ac.rs}}
\date{}

\begin{document}
\maketitle
\begin{abstract}
In this paper  some new characterizing theorems of exponential
distribution based on order statistics are presented. Some existing results are generalized  and  the open conjecture by Arnold and Villasenor is solved.
\end{abstract}

{\small \textbf{ keywords:} Characterization, order statistics, exponential distribution, Stirling numbers of the second kind

\textbf{MSC(2010):} 62E10, 62G30}

\section{Introduction}
There is an abundance of characterizations of exponential distribution and among them a considerable part is based on properties of order statistics.
Most of them could be found in \cite{ahsanullah}, \cite{arnoldHuang}, \cite{galambos} and \cite{johnson}.
Recently Arnold and Villasenor \cite{arnold} proposed a series of characterizations based on the  sample of size two and stated some conjectures on their generalization.
They also proposed a new method of proof which can be used when the density in question is analytic. Later Yanev and Chakraborty \cite{yanev}
proved by this method two characterization theorems concerning maximum of sample of size three as well as the characterization based on consecutive maxima \cite{chakraborty}.

In this paper we extend the generalizations to arbitrary order statistics. We consider the case of consecutive order statistics via convolution with
independent random variable from the same distribution. Similar problems have been studied in
\cite{weselowski} and \cite{castano}. Their formulations are slightly different in terms that the convolution in question includes a random variable with
fixed distribution.

The other case we consider is the characterization based  on representation of $k$ th order statistic of sample  of size $n$ as a weighted sum of $k$
sample members.
This problem has a history. Ahsanullah and Rahman \cite{ahsanullahRahman} proved the theorem when the representation is valid for all $k$. Later Huang \cite{huang}
showed that the condition in general cannot be relaxed to just one value of $k$.
We prove that under the assumption of analyticity of the density function the theorem is valid even under this relaxed condition.
 As a corollary we solve the  conjecture of Arnold and Villasenor regarding representation of sample maximum  stated in \cite{arnold}.

\section{Auxiliary results}
In this section we present four combinatorial identities that will be used in the proofs of characterization theorems.
In all of them appear Stirling numbers of the second kind.
A Stirling number of second kind, denoted $\stirling{a}{b}$, represents the number of ways to partition a set of $a$ objects into $b$ non-empty subsets.  In proofs of our lemmas we use the following well-known identities (see e.g. \cite{gross}).

\begin{eqnarray}\label{StirlingRekurentno}
&&\stirling{a}{b}=\stirling{a-1}{b-1}+b\stirling{a-1}{b},\\
\label{StirlingBinom}
 &&\stirling{a+1}{b+1}=\sum\limits_{l=0}^{a}\binom{a}{l}\stirling{l}{b},\\
 &&\label{StirlingDijagonalno}\stirling{a+b+1}{b}=\sum\limits_{l=0}^{b}l\stirling{a+l}{l},\\
 &&\label{StirlingStepen}a^b=\sum\limits_{l=0}^{b}\stirling{b}{l}a(a-1)\cdots(a-l+1).
 \end{eqnarray}
We proceed with the lemmas necessary for the proofs of the characterization theorems.
\begin{lemma}\label{NasaLema}
For integers $k,n,r$  such that $1< k \leq n$ and $r\geq 0$ it holds 
\begin{eqnarray}\label{NasId}
\nonumber&&\sum\limits_{j=k-2}^{k+r-1}\sum\limits_{i=0}^{j-k+2}\!\!\!\!\binom{n-k}{i} (i+k-2)!\stirling{j\!+\!1}{i\!+k\!-1}(k-1)n^{k+r-1-j}\\&=&\sum\limits_{i=0}^{r+1}\!\!\binom{n-k}{i}(i+k-1)!\stirling{k+r+1}{i+k}.
\end{eqnarray}
\end{lemma}
\textbf{Proof.} We prove the lemma by induction on $r$. For $r=0$ the equality \eqref{NasId} simplifies to
\begin{equation*}
(k-1)!\bigg(n+\stirling{k}{k-1}+(k-1)(n-k)\bigg)=(k-1)!\bigg(\stirling{k+1}{k}+k(n-k)\bigg),
\end{equation*}
which is true because of \eqref{StirlingRekurentno}.
Thus the statement of the lemma holds for $r=0$ for all  $1< k \leq n$.

Let us now suppose that \eqref{NasId} is satisfied for  $r-1$ for all  $1< k \leq n$.
 The left hand side of \eqref{NasId} can be split as
\begin{eqnarray*}
&&\nonumber\sum\limits_{j=k-2}^{k+r-2}\sum\limits_{i=0}^{j-k+2}\binom{n-k}{i} (i+k-2)!\stirling{j+1}{i+k-1}(k-1)n^{k+r-1-j}\\
&+&\sum\limits_{i=0}^{r+1}\binom{n-k}{i} (i+k-2)!\stirling{k+r}{i+k-1}(k-1).
\end{eqnarray*}
Using the induction hypothesis on the first summand we have that the expression above  is equal to
\begin{equation}\label{l2}
\sum\limits_{i=0}^{r}\binom{n-k}{i} (i+k-1)!\stirling{k+r}{i+k}n+\sum\limits_{i=0}^{r+1}\binom{n-k}{i} (i+k-2)!\stirling{k+r}{i+k-1}(k-1).
\end{equation}
It remains to prove
that
\eqref{l2} is equal to
\begin{equation*}
\sum\limits_{i=0}^{r+1}\binom{n-k}{i} (i+k-1)!\stirling{k+r+1}{i+k},
\end{equation*}
which can be written as
\begin{equation}\label{d2}
\sum\limits_{i=0}^{r+1}\binom{n-k}{i} (i+k-1)!\stirling{k+r}{i+k-1}+\sum\limits_{i=0}^{r}\binom{n-k}{i} (i+k)!\stirling{k+r}{i+k}.
\end{equation}
Grouping the corresponding summands from \eqref{l2} and \eqref{d2} we get
\begin{equation*}
\sum\limits_{i=0}^{r}\binom{n-k}{i}(i+k-1)!\stirling{k+r}{i+k}(n-k-i)=\sum\limits_{i=0}^{r+1}\binom{n-k}{i}(i+k-2)!\stirling{k+r}{i+k-1}i.
\end{equation*}
The last equality is easily shown putting $j=i+1$ in the first sum.\hfill$\Box$
\begin{lemma}\label{NasaLema2}
For integers $k,n,r$  such that $1< k \leq n$ and $r\geq 0$ it holds 
\begin{eqnarray}\label{NasId2}
&&\nonumber\sum\limits_{j=k-2}^{k+r-1}\sum\limits_{i=0}^{j-k+2}\!\!\!\!\binom{n-k+1}{i} (i+k-2)!\stirling{j\!+\!1}{i\!+k\!-1}(k-1)(n-k+1)^{k+r-1-j}\\&=&\sum\limits_{i=0}^{r+1}\!\!\binom{n-k}{i}(i+k-1)!\stirling{k+r+1}{i+k},
\end{eqnarray}
\end{lemma}
The proof of this lemma is analogous to the proof of lemma \ref{NasaLema} so we omit it here.

\begin{lemma}\label{lepaLema}
 For integers $k,n,r$  such that $1< k\leq n$ and $r\geq 0$ it holds
\begin{eqnarray}\label{lepId}
&&\nonumber \sum\limits_{i=0}^{r+1}\binom{n-k}{i}\sum\limits_{s=1}^{k}(n-k+s)\frac{(i+s-1)!}{(s-1)!}\stirling{s+r}{i+s}\\&=&\sum\limits_{i=0}^{r+1}
\binom{n-k}{i}\frac{(i+k-1)!}{(k-1)!}\stirling{k+r+1}{i+k}.
\end{eqnarray}
\end{lemma}
\textbf{Proof.}
We prove the lemma by induction on $n$. For any $r$ and $k$ and $n=k$ the expression \eqref{lepId} simplifies to identity \eqref{StirlingDijagonalno}.

Suppose now that the equality \eqref{lepId} is true for any $k$, any $r$ and $n-1$. We need to prove that it is also true  for $n$. Transforming the left hand side of \eqref{lepId} we get

\begin{eqnarray*}
&&\sum\limits_{i=0}^{r+1}\binom{n-k}{i}\bigg(\!(n-k-i)\sum\limits_{s=1}^k\frac{(i+s-1)!}{(s-1)!}\stirling{s+r}{i+s}+
\!\sum\limits_{s=1}^k\frac{(i+s)!}{(s-1)!}\stirling{s+r}{i+s}\!\bigg)
\\&=&
\sum\limits_{i=0}^{r+1}(n-k)\binom{n-k-1}{i}\sum\limits_{s=1}^k\frac{(i+s-1)!}{(s-1)!}\stirling{s+r}{i+s}\\
&+&\sum\limits_{i=0}^{r+1}\!\binom{n-k-1}{i}\!\!\sum\limits_{s=1}^k\!\!\frac{(i+s)!}{(s-1)!}\stirling{s+r}{i+s}\!+
\!\sum\limits_{i=0}^{r+1}\!\!\binom{n-k-1}{i-1}\!\sum\limits_{s=1}^k\!\frac{(i+s)!}{(s-1)!}\stirling{s+r}{i+s}\\
&=&\sum\limits_{i=0}^{r+1}(n-k)\binom{n-k-1}{i}\sum\limits_{s=1}^k\frac{(i+s-1)!}{(s-1)!}\stirling{s+r}{i+s}\\
&+&\sum\limits_{i=0}^{r+1}\binom{n-k-1}{i}\sum\limits_{s=1}^k\frac{(i+s)!}{(s-1)!}\bigg(\stirling{s+r}{i+s}+(i+s+1)\stirling{s+r}{i+1+s}\bigg).
\end{eqnarray*}
Applying the identity \eqref{StirlingRekurentno}, shifting the index $s$ in the last inner sum and separating  the term for $s=k+1$, the expression above 
becomes

\begin{eqnarray*}
&&\sum\limits_{i=0}^{r+1}(n-k)\binom{n-k-1}{i}\sum\limits_{s=1}^k\frac{(i+s-1)!}{(s-1)!}\stirling{s+r}{i+s}\\&+&
\sum\limits_{i=0}^{r+1}\binom{n-k-1}{i}\sum\limits_{s=2}^k\frac{(i+s-1)!}{(s-2)!}\stirling{s+r}{i+s}\\
&+&\sum\limits_{i=0}^{r+1}\binom{n-k-1}{i}\frac{(i+k)!}{(k-1)!}\stirling{k+1+r}{i+k+1}.
\end{eqnarray*}
Grouping the first two summands together and applying the induction hypothesis to the result  we get the right hand side
of \eqref{lepId}. 
\hfill{$\Box$}
\begin{lemma}\label{ruznaLema}
 For integers $k,n,r$  such that  $1<k\leq n$ and $r\geq 0$ it holds
\begin{equation}\label{ruzniId}
\sum\limits_{\underset{j_1+\cdots+j_{k}=r+1}{j_1,\ldots,j_{k}\geq 0}}n^{j_1}(n-1)^{j_2}\cdots(n-k+1)^{j_k}=
\sum\limits_{i=0}^{r+1}\binom{n-k}{i}\frac{(i+k-1)!}{(k-1)!}\stirling{i+r+1}{i+k}.
\end{equation}
\end{lemma}
\textbf{Proof.} The proof is done using the strong induction on $r$. For any $k$ and $n$ and $r=0$ we have
\begin{equation*}
n+(n-1)+\cdots+(n-k+1)=\stirling{k+1}{k}+(n-k)k,
\end{equation*}
which is obviously true since $\stirling{k+1}{k}=\frac{k(k+1)}{2}$. Suppose now that \eqref{ruzniId}
is satisfied up to $r-1$. Then it remains to prove that it is satisfied for $r$.

Splitting the sum on the left hand side of \eqref{ruzniId} into two parts: for $j_1=0$ and $j_1\geq 0$, we get
\begin{eqnarray}\label{brake}
\!\!\!\!\!&&\nonumber\sum\limits_{\underset{j_1+\cdots+j_{k}=r+1}{j_1,\ldots,j_{k}\geq 0}}\!\!\!n^{j_1}(n-1)^{j_2}\cdots(n-k+1)^{j_k}=
\!\!\!\!\!\sum\limits_{\underset{j_1+\cdots+j_{k}=r+1}{j_1,\ldots,j_{k}\geq 0}}\!\!\!(n-1)^{j_2}\cdots(n-k+1)^{j_k}
\\\!\!\!\!\!&&+\sum\limits_{j_1=1}^{r+1}n^{j_1}\sum\limits_{\underset{j_2+\cdots+j_{k}=r+1-j_1}{j_2,\ldots,j_{k}\geq 0}}(n-1)^{j_2}\cdots(n-k+1)^{j_k}.
\end{eqnarray}
The sum of indices in the inner sum of the second summand of \eqref{brake} is smaller than $r+1$ so the induction hypothesis is
applicable (in this case for $n-1$, $k-1$ and $r+1-j_1$). The first summand can be recursively split in the same manner until all indices except the last one are  equal to zero. After this process (including the application of induction hypothesis) we obtain

\begin{eqnarray}\label{jpresmene}
 \!\!\!\!\!\!\!&& \sum\limits_{\underset{j_1+\cdots+j_{k}=r+1}{j_1,\ldots,j_{k}\geq 0}}n^{j_1}(n-1)^{j_2}\cdots(n-k+1)^{j_k}=\;(n-k+1)^{r+1}+\\\nonumber\!\!\!\!\!\!\!&&+
\sum\limits_{l=1}^{k-1}\sum\limits_{j=1}^{r+1}(n-l+1)^{j}\sum\limits_{i=0}^{r+1-j}\binom{n-k}{i}\frac{(i+k-l-1)!}{(k-l-1)!}\stirling{k-l+r+1-j}{i+k-l}.
\end{eqnarray}
Substituting the index $j$ with $m=k+r-1-l-j$ and, subsequently, the index $l$ with $s=k-l+1$, as well as applying the identity \eqref{StirlingStepen} to $(n-k+1)^{r+1}$, \eqref{jpresmene} becomes
\begin{eqnarray*}
&&\sum\limits_{s=2}^{k}\sum\limits_{i=0}^{r}\sum\limits_{m=i+s-2}^{r+s-2}\binom{n-k}{i}(n-k+s)^{r+s-1-m}\frac{(i+s-2)!}{(s-2)!}
\stirling{m+1}{i+s-1}\\&+&\sum\limits_{i=0}^{r+1}\frac{(n-k+1)!}{(n-k+1-i)!}\stirling{r+1}{i}
\\&=&\sum\limits_{s=2}^k\frac{(n-k+s)}{(s-1)!}\sum\limits_{i=0}^{r}\binom{n-k}{i}\sum\limits_{m=i+s-2}^{r+s-2}\bigg(n-k+s)^{r+s-2-m}(s-1)\\
&&(i+s-2)!
\stirling{m+1}{i+s-1}\bigg)+(n-k+1)\sum\limits_{i=1}^{r+1}\binom{n-k}{i-1}(i-1)!\stirling{r+1}{i}.
\end{eqnarray*}
Applying lemma \ref{NasaLema} to the two inner sums and grouping the summands we get
\begin{equation*}
\!\!\!\sum\limits_{s=1}^{k}\sum\limits_{i=0}^{r}(n-k+s)\binom{n-k}{i}\frac{(i+s-1)!}{(s-1)!}\stirling{r+s}{i+s}.
\end{equation*}
Applying now lemma \ref{lepaLema} we obtain the right hand side of \eqref{ruzniId}. Hence  the proof is completed.\hfill{$\Box$}
\begin{remark}
 The statements of lemmas \ref{lepaLema} and \ref{ruznaLema} are also true for $k=1$. They  could be easily proven analogously. 
\end{remark}

\section{Main results}
In the beginning we state and prove two lemmas that will play an  important role in the proofs of the theorems. They are similar to those from \cite{chakraborty}.

Let $\mathcal{F}$ be a class of continuous distribution functions $F$ such that $F(0)=0$ and whose density function $f$ allows expansion in Maclaurin series for all $x>0$.

\begin{lemma}\label{lemaOsnova}
 Let $F$ be a distribution function that belongs to $\mathcal{F}$. If for all natural $q$ holds
\begin{equation}\label{fizvod}
  f^{(q)}(0)=(-1)^qf^{q+1}(0),
 \end{equation}
 then $f(x)=\lambda e^{-\lambda x}$ for some $\lambda>0$.
\end{lemma}
\textbf{Proof.} Expanding the function $f$ in Maclaurin series for positive values of $x$ we get
\begin{equation}
 f(x)=\sum\limits_{q=0}^{\infty}f^{(q)}(0)\frac{x^{q}}{q!}=\sum\limits_{q=0}^{\infty}(-1)^qf^{q+1}(0)\frac{x^{q}}{q!}=f(0)e^{-f(0)x}.
\end{equation}
For $f(0)>0$ this is the density of exponential distribution with $\lambda=f(0)$. $\hfill \Box$

\begin{lemma}\label{lemaAizvodi} Let $F$ be a distribution function that belong to the class $\mathcal{F}$. Denote $A_{m}(x) = F^{m}(x)f(x)$.
If the condition \eqref{fizvod} is satisfied for all $0\leq k\leq r-m$, $r>m$, then
\begin{equation}\label{Amizvod}
 A_m^{(r)}(0)=(-1)^{r-m}f^{r+1}(0)\stirling{r+1}{m+1}m!.
\end{equation}
\end{lemma}
\begin{remark}
 For $r\leq m$ the statement of the lemma is valid without any condition imposed on derivatives of $f$.
\end{remark}
{\bf Proof.} The $r$ th derivative of $A_m(x)$ is

\begin{equation*}
 A_m^{(r)}(x)=\sum\limits_{\underset{j_1+\cdots+j_{m+1}=r}{j_1,\ldots,j_{m+1}\geq 0}}\binom{r}{j_1,\ldots,j_{m+1}}F^{(j_1)}(x)\cdots F^{(j_m)}(x) f^{(j_{m+1})}(x).
\end{equation*}
Using the fact that $F(0)=0$ we get

\begin{equation}\label{idAm}
 A_m^{(r)}(0)=\sum\limits_{\underset{j_1+\cdots+j_{m+1}=r}{j_1,\ldots,j_m\geq 1,\;j_{m+1}\geq 0}}\binom{r}{j_1,\ldots,j_{m+1}}f^{(j_1-1)}(0)\cdots f^{(j_m-1)}(0) f^{(j_{m+1})}(0).\\
 \end{equation}
 Since all derivatives are of orders smaller or equal to $r-m$ using \eqref{fizvod} we obtain
\begin{eqnarray*} \nonumber A_m^{(r)}(0) &=&\sum\limits_{\underset{j_1+\cdots+j_{m+1}=r}{j_1,\ldots,j_m\geq 1,\;j_{m+1}\geq 0}}\binom{r}{j_1,\ldots,j_{m+1}}(-1)^{r-m}f^{r+1}(0)\\
             \nonumber&=&(-1)^{r-m}f^{r+1}(0)\sum\limits_{\underset{j_1+\cdots+j_{m+1}=r}{j_1,\ldots,j_m,j_{m+1}\geq 1}}\binom{r}{j_1,\ldots,j_{m+1}}\\
            \nonumber &+& (-1)^{r-m}f^{r+1}(0)\sum\limits_{\underset{j_1+\cdots+j_{m+1}=r}{j_1,\ldots,j_m\geq 1,\;j_{m+1}= 0}}\binom{r}{j_1,\ldots,j_m}\\
             \nonumber&=& (-1)^{r-m}f^{r+1}(0)\bigg(\stirling{r}{m+1}(m+1)!+\stirling{r}{m}m!\bigg)\\
             \nonumber&=& (-1)^{r-m}f^{r+1}(0) m! \stirling{r+1}{m+1}.
\end{eqnarray*}
In the last line we used the  identity \eqref{StirlingRekurentno}.\hfill{$\Box$}
\bigskip

Let $X_{(k;n)}$ be the $k$ th order statistics from the sample of size $n$.
We now state the characterization theorems.
\begin{theorem}\label{nasaGlavna}
Let $X_1,\ldots,X_n$ be a random sample from the distribution $F$ that belongs to $\mathcal{F}$. Let $k$ be a fixed number
such that $1<k\leq n$. If

\begin{equation}\label{glavna}
 X_{(k-1;n-1)}+\frac{1}{n}X_n\overset{d}{=} X_{(k;n)}
\end{equation}
then $X\sim \mathcal{E}(\lambda), \lambda>0$.
\end{theorem}
{\bf Proof. }
Equalizing the  densities from \eqref{glavna} we get
\begin{eqnarray*}
&&\int\limits_{0}^{x}\frac{(n-1)!}{(k-2)!(n-k)!}F^{k-2}(x-y)(1-F(x-y))^{n-k}f(x-y)nf(ny)dy\\&&=\frac{n!}{(k-1)!(n-k)!}F^{k-1}(x)(1-F(x))^{n-k}f(x),
\end{eqnarray*}
or

\begin{eqnarray}\label{integralnaJ}
\nonumber&&(k-1)\sum\limits_{i=0}^{n-k}(-1)^{i}\binom{n-k}{i}\int\limits_{0}^{x}A_{i+k-2}(x-y)f(ny)dy\\&=&f(x)\sum\limits_{i=0}^{n-k}(-1)^{i}\binom{n-k}{i}\int\limits_{0}^{x}A_{i+k-2}(y)dy.
\end{eqnarray}
Using induction we prove that \eqref{fizvod} holds for every natural $q$ which by lemma \ref{lemaOsnova} implies that $f(x)$ is exponential density.






Differentiating integral equation \eqref{integralnaJ} $k$ times we get

\begin{eqnarray*}
&&(k-1)\sum\limits_{i=0}^{n-k}(-1)^{i}\binom{n-k}{i}
\bigg(\sum\limits_{j=0}^{k-1}n^{k-1-j}f^{(k-1-j)}(nx)A^{(j)}_{i+k-2}(0)
\\&+&\int\limits_{0}^{x}A^{(k)}_{i+k-2}(x-y)f(ny)dy\bigg)\\
&=&\sum\limits_{i=0}^{n-k}(-1)^{i}(i+k-1)\binom{n-k}{i}
\bigg(\sum\limits_{j=1}^{k}\binom{k}{j}f^{(k-j)}(x)A^{(j-1)}_{i+k-2}(0)
\\&+&\!f^{(k)}(x)\int\limits_{0}^{x}A^{(k)}_{i+k-2}(y)dy\bigg).
\end{eqnarray*}

Letting $x=0$  and eliminating zero terms we get

\begin{eqnarray*}
 \!\!\!\!\!\!\!\!\!&&\!(k-1)\big(nf'(0)(k-2)!f^{k-1}(0)\!+\!f(0)A^{(k-1)}_{k-2}(0)\!-\!(n-k)(k-1)!f^{k+1}(0)\big)=\\
 \!\!\!\!\!\!\!\!\!&&\!(k-1)f(0)A^{(k-1)}_{k-2}(0)\!+\!(k-1)f'(0)(k-2)!f^{k-1}(0)k\!-\!k(n-k)(k-1)!f^{k+1}(0),
\end{eqnarray*}
from where we get $f'(0)=-f^2(0),$
which means that \eqref{fizvod} holds for $q=1$. Suppose now that \eqref{fizvod} is satisfied for all $k\leq r$. We shall prove that it holds for $q=r+1$.

Differentiating the  integral equation \eqref{integralnaJ} $k+r$ times, letting $x=0$ and eliminating zero terms we get 
\begin{eqnarray*}
 &&(k-1)\sum\limits_{i=0}^{r+1}(-1)^{i}\binom{n-k}{i}\sum\limits_{j=i+k-2}^{k+r-1}n^{k+r-1-j}f^{(k+r-1-j)}(0)A^{(j)}_{i+k-2}(0)
 \\&&=\sum\limits_{i=0}^{r+1}(-1)^{i}(i+k-1)\binom{n-k}{i}\sum\limits_{j=i+k-2}^{k+r-1}\binom{k+r}{j+1}f^{(k+r-1-j)}(0)A^{(j)}_{i+k-2}(0).
\end{eqnarray*}
The terms for $i=0$ and $j=k+r-1$ are equal and hence they cancel out. Splitting the summation into two parts  for $i=0$ and $i>0$ we get

\begin{eqnarray*}
\!\!\!\!\!\!\!\! \!\! &&\!\!(k-1)\bigg(n^{r+1}f^{(r+1)}(0)A^{(k-2)}_{k-2}(0)+\sum\limits_{j=k-1}^{k+r-2}n^{k+r-1-j}f^{(k+r-1-j)}(0)A^{(j)}_{k-2}\bigg)
\!\!\!\!\!\!\!\!\\&+&\!\!(k-1)\sum\limits_{i=1}^{r+1}\sum\limits_{j=i+k-2}^{k+r-1}(-1)^{i}\binom{n-k}{i}n^{k+r-1-j}f^{(k+r-1-j)}(0)A^{(j)}_{i+k-2}(0)
\!\!\!\!\!\!\!\!\\&=&\!\!(k-1)\binom{k+r}{k-1}f^{(r+1)}(0)A^{(k-2)}_{k-2}(0)
\!\!\!\!\!\!\!\!\\&+&\!\!(k-1)\sum\limits_{k-1}^{k+r-2}\binom{k+r}{j+1}f^{(k+r-1-j)}(0)A^{(j)}_{k-2}(0)
\!\!\!\!\!\!\!\!\\&+&\!\!\sum\limits_{i=1}^{r+1}\sum\limits_{j=i+k-2}^{k+r-1}(-1)^{i}\binom{n-k}{i}(i+k-1)\binom{k+r}{j+1}f^{(k+r-1-j)}(0)A^{(j)}_{i+k-2}(0).
\end{eqnarray*}
Applying the induction hypothesis to the derivatives of functions $f$ and, consequently, via lemma \ref{lemaAizvodi}, to  $A_{i+k-2}$ and grouping the summands we obtain
\begin{eqnarray*}
 &&f^{(r+1)}(0)f^{k-1}(0)(k-1)!\bigg(n^{r+1}-\binom{k+r}{k-1}\bigg)\\
 &=&(-1)^{r+1}f^{k+r+1}(0)\Bigg((k-1)!\sum\limits_{j=k-1}^{k-2}\Big(\binom{k+r}{j+1}-n^{k+r-1-j}\Big)\stirling{j+1}{k-1}
 \\&+&\sum\limits_{i=1}^{r+1}\sum\limits_{j=i+k-2}^{k+r-1}\binom{n-k}{i}\Big((i+k-1)\binom{k+r}{j+1}\\&-&
 (k-1)n^{k+r-1-j}\Big)(i+k-2)!\stirling{j+1}{i+k-1}\Bigg).
\end{eqnarray*}
To prove the induction step it remains to show that
\begin{eqnarray*}
 &&(k-1)!\bigg(n^{r+1}-\binom{k+r}{k-1}+\sum\limits_{j=k-1}^{k-2}\Big(n^{k+r-1-j}-\binom{k+r}{j+1}\Big)\stirling{j+1}{k-1}\bigg)\\
 &=&\sum\limits_{i=1}^{r+1}\sum\limits_{j=i+k-2}^{k+r-1}\binom{n-k}{i}\Big((i+k-1)\binom{k+r}{j+1}\\&-&
 (k-1)n^{k+r-1-j}\Big)(i+k-2)!\stirling{j+1}{i+k-1}.
\end{eqnarray*}

Joining the summation for $i=0$ and $i>0$ back together we get
\begin{eqnarray*}
&&\sum\limits_{i=0}^{r+1}\sum\limits_{j=i+k-2}^{k+r-1}\binom{n-k}{i}(k-1)n^{k+r-1-j}(i+k-2)!\stirling{j+1}{i+k-1}\\
&&=\sum\limits_{i=0}^{r+1}\binom{n-k}{i}(i+k-1)!\sum\limits_{j=i+k-2}^{k+r-1}\binom{k+r}{j+1}\stirling{j+1}{i+k-1}.
\end{eqnarray*}

Using  identity  \eqref{StirlingBinom} and lemma \ref{NasaLema} we complete the proof.\hfill{$\Box$}
\begin{theorem}\label{nasaGlavna2}
Let $X_1,\ldots,X_n$ be a random sample from the distribution $F$ that belongs to $\mathcal{F}$ and let $X_0$ be a random variable  independent of the sample that follows the same distribution. Let $k$ be a fixed number
such that $1<k\leq n$. If

\begin{equation}\label{glavna2}
 X_{(k-1;n)}+\frac{1}{n-k+1}X_0\overset{d}{=} X_{(k;n)}
\end{equation}
then $X\sim \mathcal{E}(\lambda), \lambda>0$.
\end{theorem}
We omit the proof since it follows  completely analogous  procedure to the proof of theorem \ref{nasaGlavna} with the application of lemma \ref{NasaLema2} in the last step.

\begin{theorem}\label{nasaGlavna3}
Let $X_1,\ldots,X_n$ be a random sample from the distribution $F$ that belongs to $\mathcal{F}$. Let $k$ be a fixed number
such that $1\leq k\leq n$. If

\begin{equation}\label{glavna3}
 \frac{1}{n}X_1+\frac{1}{n-1}X_2+\cdots+\frac{1}{n-k+1}X_k\overset{d}{=} X_{(k;n)}
\end{equation}
then $X\sim \mathcal{E}(\lambda), \lambda>0$.
\end{theorem}
{\bf Proof. } Let $k\geq 2$.
Equalizing  the respective densities as in the previous proof  we get

\begin{eqnarray}\label{SanIdentitet}
\nonumber &&\int\limits_{0}^x\!f(n(x-y_2))\cdots\!\int\limits_{0}^{y_{k-1}}\!\!\!f((n-k+2)(y_{k-1}-y_{k}))f((n-k+1)y_k)dy_2\cdots dy_k
\\&&=\frac{1}{(k-1)!}f(x)\sum\limits_{i=0}^{n-k}(-1)^{i}\binom{n-k}{i}\int\limits_{0}^{\infty}A_{k-2+i}(x)dx.
\end{eqnarray}
Denote  the left hand side of \eqref{SanIdentitet} with $J_{k,n}(x)$. Obviously, it can be expressed as
\begin{eqnarray*}
&&J_{k,n}(x)=\int\limits_{0}^xf(n(x-y_2))J_{k-1,n-1}(y_2)dy_2,\\
&&J_{1,1}(x)=f((n-k+1)x).
\end{eqnarray*}
The $(k+r)$ th derivative of $J_{k,n}$ is
\begin{eqnarray*}
J_{k,n}^{(k+r)}(x)&=&\sum\limits_{j=0}^{k+r-1}n^jf^{(j)}(0)J^{(k+r-j-1)}_{k-1,n-1}(x)
\\&+&\int\limits_{0}^x
f^{(k+1)}(n(x-y_2))n^{k+r}J^{(r+1)}_{k-1,n-1}(y_2)dy_2.\end{eqnarray*}
Letting $x=0$ we get
\begin{eqnarray}
\label{jknr}
\nonumber&&J_{k,n}^{(k+r)}(0)=\sum\limits_{j=0}^{k+r-1}n^jf^{(j)}(0)J^{(k+r-j-1)}_{k-1,n-1}(0),\\&&J^{(s)}_{1,1}(0)=(n-k+1)^sf^{(s)}(0),\text{ for every $s\geq 0$.}
\end{eqnarray}
Applying the recurrence relation \eqref{jknr} $k-1$ times we obtain

\begin{eqnarray*}
&&J^{(k+1)}_{k,n}(0)=\sum\limits_{j_1=0}^{k+r-1}n^{j_1}f^{(j_1)}(0)\sum\limits_{j_2=0}^{k+r-2-j_1}(n-1)^{j_2}f^{(j_2)}(0)\cdots\\&\cdots&\!\!\!
\sum\limits_{j_{k-1}=0}^{r+1-
\sum\limits_{l=1}^{k-2}j_l}(n-k+2)^{j_{k-1}}f^{(j_{k-1})}(n-k+1)^{r+1-\sum\limits_{l=1}^{k-1}j_l}f^{(r+1-\sum\limits_{l=1}^{k-1}j_l)}(0).
\end{eqnarray*}
Then the $(k+r)$ th derivative of the left hand side of \eqref{SanIdentitet} becomes
\begin{equation*}
\sum\limits_{\underset{j_1+\cdots+j_{k}=r+1}{j_1,\ldots,j_{k}\geq 0}}n^{j_1}(n-1)^{j_2}\cdots(n-k+1)^{j_k}
f^{(j_1)}(0)f^{(j_2)}(0)\cdots f^{(j_k)}(0).
\end{equation*}

As before we shall  prove by induction that \eqref{fizvod} holds for every $q$.
For $r=0$, the $k$ th derivative   of \eqref{SanIdentitet} at $x=0$ is
\begin{eqnarray}\label{ind1}
\nonumber&& (n+n-1+\cdots+n-k+1)f'(0)f^{k-1}(0)\\&=&\frac{1}{(k-2)!}f(0)A^{(k-1)}_{k-2}(0)+f'(0)f^{k-1}(0)k-k(n-k)f^{k+1}(0).
\end{eqnarray}
From \eqref{idAm} we can get that
\begin{equation*}
A^{(k-1)}_{k-2}(0)=f^{k-2}(0)f'(0)(k-1)!+(k-2)f'(0)f^{k-2}(0)\frac{(k-1)!}{2}.
\end{equation*}
Inserting this in \eqref{ind1} we get $f'(0)=-f^{2}(0)$ which means \eqref{fizvod} holds for $q=1$.
Suppose now that \eqref{fizvod} is satisfied for all $q\leq r$. We shall prove that it holds for $q=r+1$.
The $(k+r)$ th derivative  of \eqref{SanIdentitet} at $x=0$ is

\begin{eqnarray}\label{idizvodi}
\!\!&&\!\!\!\sum\limits_{\underset{j_1+\cdots+j_{k}=r+1}{j_1,\ldots,j_{k}\geq 0}}n^{j_1}(n-1)^{j_2}\cdots(n-k+1)^{j_k}
f^{(j_1)}(0)f^{(j_2)}(0)\cdots f^{(j_k)}(0)
\\\nonumber\!\!&=&\!\!\!\sum\limits_{i=0}^{r+1}(-1)^i\frac{(i+k-1)}{(k-1)!}
\binom{n-k}{i}\sum\limits_{j=i+k-2}^{k+r-1}\binom{k+r}{j+1}f^{(k+r-1-j)}(0)A^{(j)}_{i+k-2}(0).
\end{eqnarray}
Applying the induction hypothesis the left hand side of  \eqref{idizvodi} becomes
\begin{eqnarray*}
&&f^{k-1}(0)f^{(r+1)}(n^{r+1}+\cdots+(n-k+1)^{r+1})\\&+&\sum\limits_{\underset{j_1+\cdots+j_{k}=r+1}{0\leq j_1,\ldots,j_{k}<r+1 }}n^{j_1}(n-1)^{j_2}\cdots(n-k+1)^{j_k}(-1)^{r+1}f^{r+1+k}(0),
\end{eqnarray*}
while the right hand side of  \eqref{idizvodi} can be expressed as

\begin{eqnarray*}
\!\!&&\!\!\!\!\sum\limits_{i=1}^{r+1}\binom{n-k}{i}\sum\limits_{j=i+k-2}^{k+r-1}\binom{k+r}{j+1}(-1)^{r+1}f^{k+r+1}\frac{(i+k-1)!}{(k-1)!}\stirling{j+1}{i+k-1}
\\\!\!&+&\!\!\!\!\sum\limits_{j=i+k-1}^{k+r-2}f^{k+r+1}(0)\frac{(i+k-2)!}{(k-2)!}(-1)^{r+1}\stirling{j+1}{k-1}
\\&+&\binom{k+r}{k-1}f^{(r+1)}(0)\frac{(i+k-2)!}{(k-2)!}
+\frac{1}{(k-2)!}f(0)A^{(k+r-1)}_{k-2}(0).
\end{eqnarray*}
The term  $A^{(k+r-1)}_{k-2}(0)$ can be evaluated using \eqref{Amizvod} and \eqref{idAm}  as
\begin{eqnarray*}
A^{k+r-1}_{k-2}(0)&=&\frac{(k+r-1)!}{(r+1)!}f^{k-2}(0)f^{(r+1)}(0)
\\&+&\frac{(k+r-1)!}{(r+2)!}(k-2)f^{k-2}(0)f^{(r+1)}(0)\\
&+&\sum\limits_{\underset{j_1+\cdots+j_{k}=k+r-1}{\underset{0\leq j_{k-1}<r+1}{1\leq j_1,\ldots,j_{k-2}<r+2}}}
(-1)^{r+1}f^{r+k}(0)\frac{(k+r-1)!}{j_1!\cdots j_{k-1}!}.
\end{eqnarray*}
After transformations given above and grouping the summands \eqref{idizvodi} becomes
\begin{eqnarray*}
\!\!&&\!\!\!\!\!\!f^{(r+1)}(0)\Bigg(\!n^{r+1}\!+\!\cdots\!+\!(n-k+1)^{r+1}\!-\!\binom{k+r}{k-1}\!-\!\binom{k+r-1}{k-2}\!-\!\binom{k+r-1}{k-3}\!\Bigg)
\\\!\!&=&\!\!\!\!\!(-1)^{r+1}f^{r+2}(0)\Bigg(\sum\limits_{j=k-1}^{k+r-2}\binom{k+r}{j+1}\stirling{j+1}{k-1}
\\\!\!&+&\!\!\!\!\!\!\sum\limits_{i=1}^{r+1}\binom{n-k}{i}\sum\limits_{j=i+k-2}^{k+r-1}\binom{k+r}{j+1}\frac{(i+k-1)!}{(k-1)!}\stirling{j+1}{i+k-1}
\\\!\!&+&\!\!\!\!\!\!\frac{1}{(k-2)!}\!\!\sum\limits_{\underset{j_1+\cdots+j_{k}=k+r-1}{\underset{0\leq j_{k-1}< r+1}{1\leq j_1,\ldots,j_{k-2}<r+2}}}\frac{(k+r-1)!}{j_1!\cdots j_{k-1}!}
\!-\!\!\!\sum\limits_{\underset{j_1+\cdots+j_{k}=r+1}{0\leq j_1,\ldots,j_{k}<r+1}}\!\!\!\!\!\!\!\!\!n^{j_1}(n-1)^{j_2}\cdots(n-k+1)^{j_k}\Bigg).
\end{eqnarray*}
To prove the induction step it remains to show that the expressions in the brackets on both sides are equal.
Joining the summands back together and applying the identity \eqref{StirlingBinom} we obtain

\begin{equation*}
\sum\limits_{\underset{j_1+\cdots+j_{k}=r+1}{j_1,\ldots,j_{k}\geq 0}}n^{j_1}(n-1)^{j_2}\cdots(n-k+1)^{j_k}=
\sum\limits_{i=0}^{r+1}\binom{n-k}{i}\frac{(i+k-1)!}{(k-1)!}\stirling{i+r+1}{i+k},
\end{equation*}
which follows from lemma \ref{ruznaLema}. Hence, the proof for $k\geq 2$ is completed.

In case of $k=1$ the proof is done in an analogous
way, but it is much simpler, so we omit it here. \hfill$\Box$
\bigskip

The following corollary, which follows directly from theorem \ref{nasaGlavna3},  is a conjecture stated in \cite{arnold}. 
\begin{corollary}
 Let $X_1,\ldots,X_n$ be a random sample from the distribution $F$ that belongs to $\mathcal{F}$. If

\begin{equation*}
 X_1+\frac{1}{2}X_2+\cdots+\frac{1}{n}X_n\overset{d}{=} X_{(n;n)},
\end{equation*}then $X\sim \mathcal{E}(\lambda), \lambda>0$.
\end{corollary}


\begin{thebibliography}{4}


\bibitem{ahsanullah} M. Ahsanullah, G.G. Hamedani, Exponential distribution: Theory and Methods, NOVA Science, New York, 2010.

\bibitem{ahsanullahRahman} M. Ahsanullah, M. Rahman, A Characterization of the Exponential Distribution,
J Appl Prob, \textbf{9}(2) (1972), 457--461.

\bibitem{arnoldHuang} B.C. Arnold, J.S. Huang, Chapter 12: Characterizations, In: N. Balakrishnan, P.A. Basu, The Exponential Distribution:
Theory, Methods and Applications, Gordon and Breach, Amsterdam, 1995. (185--203)

\bibitem{arnold} B.C. Arnold, J.A. Villasenor,  Exponential characterizations motivated by the
structure of order statistics in samples of size two, Stat Probab Lett  \textbf{83}(2) (2013), 596 -- 601.

\bibitem{castano} A. Casta\~ no-Martinez, F. L\' opez-Bl\' azquez, B. Salamanca-Mi\~no, Random translations, contractions and dilations
of order statistics and records, Statistics \textbf{46}(1) (2012), 57--67.

\bibitem{chakraborty} S. Chakraborty, G.P. Yanev, Characterization of exponential distribution through
equidistribution conditions for consecutive maxima, J Stat Appl Pro \textbf{2}(3) (2013), 237 -- 242.

\bibitem{galambos} J. Galambos, S. Kotz, Characterizations of Probability Distributions,
Springer-Verlag, Berlin-Heidelberg-New York, 1978.

\bibitem{gross} J. Gross, Combinatorial Methods with Computer Applications, Chapman \& Hall, 2008.

\bibitem{huang} J.S. Huang, On a Theorem of Ahsanullah and Rahman, J Appl Prob, \textbf{11}(1) (1974), 216--218.

\bibitem{johnson} N.L. Johnson, S. Kotz, N. Balakrishnan, Continuous Univariate Distributions, vol.1, 2nd Ed., Wiley, New York, 1994.


\bibitem{weselowski}J. Weselowski, M. Ahsanullah, Switching order statistics through random power contractions, Aust N Z J Stat,
\textbf{46}(2) (2004), 297--303.

\bibitem{yanev} G.P. Yanev, S. Chakraborty, Characterizations of exponential distribution based on sample of size three,
Pliska Stud Math Bulgar \textbf{23} (2013), 237 -- 244.

\end{thebibliography}
\end{document}